\documentclass[12pt,oneside]{article}
\usepackage{amsmath,amssymb,amsfonts,amsthm}
\def\goth{\mathfrak}
\def\Sc{\mathcal}

\def\:{{\em\,:}}
\def\({{\em (}}
\def\){{\em )}}
\def\[{{\em [}}
\def\]{{\em ]}}
\def\1#1{\big#1}
\def\2#1{\Big#1}
\def\3#1{\bigg#1}
\def\4#1{\Bigg#1}

\def\o#1{\overline{#1}}

\def\G{\Gamma}
\def\O{\Omega}
\def\x{{\goth X}(T_2M)}
\def\xo{{\goth X}^{\pi_1}(T_2M)}
\def\xt{{\goth X}^{\pi_2}(T_2M)}

\def\ttm{{\Sc T}_2M}

\def\tm{T_2M}
\def\vo{V^{\pi_1}(T_2M)}
\def\vt{V^{\pi_2}(T_2M)}

\def\Section#1{\vspace{30truept}\addtocounter{section}{1}\centerline{\Large\bf
	\arabic{section}.~~#1}\vspace{12pt}}
\newtheorem{definition}{Definition}
\newtheorem{theorem}{Theorem}
\newtheorem{lemma}{Lemma}
\newtheorem{proposition}{Proposition}
\newtheorem{corollary}{Corollary}

\newtheorem{remark}{Remark}   

\title{SPRAYS AND CONNECTIONS ON THE TANGENT BUNDLE OF ORDER TWO}
\author{Nabil L. Youssef} 
\date{}

\begin{document}               
\bibliographystyle{plain}
\maketitle                     
\vspace{-1.15cm}
\begin{center}
{Department of Mathematics, Faculty of Science,\\ Cairo
University, Giza, Egypt.} 
\end{center}
\vspace{1cm}
\maketitle                         
\smallskip

\noindent{\bf Abstract.}
Adopting the global approach to tangent bundles of order two established in [1], we develop 
this approach to find new results. We also generalize various results of [3], [4] and [6] to the
geometry of tangent bundles of order two.
\bigskip


\vspace{30truept}\centerline{\Large\bf{Introduction}}\vspace{12pt}
\par Many authors have studied the geometry of tangent bundles of order two ([1], [9], [10])
and even of higher order ([5],[7], [8]). However, global approachs to second order geometries
seem to be rare in the literature. An interesting global approach to these geometries is due to
G. CATZ [1]. In this approach she has developed a global theory of the geometry of tangent
bundles of order two $\tm$, inspired by Grifone`s theory ([3], [4])of the geometry of tangent 
bundles $TM$, but based essentially on two canonical tensor fields $J_1$ and $J_2$ of type 
$(1,1)$ on $\tm$, instead of the natural almost tangent structure $J$ on $TM$ considered by
Grifone and other authors.

It is this global point of view of Catz which we shall adopt in this paper. Our aim is in fact
twofold. On one hand, we continue the development of Catz`s approach and find new results.
On the other hand, we generalize various results of [3], [4] and [6] to the geometry of tangent
bundles of order two.

The paper consists of thre parts preceded by an introductory section (§1), which provides a 
brief account on the geometry of tangent bundles of order two as established by Catz [1].

In the first part (§2), we prove the so-called theorem of canonical decomposition of connections 
of type $1$ on $\tm$.

The second part (§3) is devoted to the study of nonlinear connections on $\tm$ in relation 
with linear connections. We show that a linear connection on $\tm$ satisfying certain 
conditions induces nonlinear connections of type $1$ and of type $2$ on $\tm$, and we 
investigate these nonlinear connections. This generalizes some results of [4] to $\tm$.

In the third and last part (§4), we define and investigate Finslerian forms $\Omega$ on $\tm$.
We study the canonical tensors associated to $\Omega$ and we establish, in paritcular, the 
existence and uniqueness of a canonical spray of type $2$, a canonical connection of type 
$2$ and a canonical connection of type $1$ on $\tm$. This generalizes to $\tm$ various 
interesting results of [3] and[6].

We point out here and once that all geometric objects considered in this paper are supposed 
to be of class $C^{\infty}$. The formalism of Fr\"{o}licher-Nijenhuis [2] will be our 
fundamental tool.

\Section{Tangent Bundle of Order Two }
\par In this section we give a brief account on the geometry of tangent bundles of order 2 as 
established by Catz [1]; we keep the same notations used in this reference. We emphasise
the definitions and the main results which will be used in the sequel.          

Let $M$ be a differentiable manifold of dimension $n$. Let $\tm$ be the set of $2$-jets with
origin $O$ determined by the differentiable mappings of $\Bbb R$ into $M$. An element of
$\tm$ is written $X=J^2_0f$, with $f:\Bbb R\longrightarrow M$.

We have the following projections:
$$\pi_1:\tm\longrightarrow TM:X=J^2_0f\longmapsto\pi_1(X)=J^1_0f$$
$$\pi_2:\tm\longrightarrow TM:X=J^2_0f\longmapsto\pi_2(X)=f(0)$$

Let $X=J^2_0f\in\tm$. If $(U,(x_1,\ldots,x_n))$ is a local chart of $M$ at $f(0)$, the jet $X$ is 
determined by
$$x_i(X)=x_i(f(0)),\quad y_i(X)=\frac{dx_i\circ f}{dt}(0),\quad z_i(X)=\frac{d^2x_i\circ f}{d^2t}(0),
\quad i=1\ldots,n.$$
The set $\tm$ equipped with the chart $(\pi^{-1}_2(U),(x_i,y_i,z_i))$ is a differentiable manifold 
of dimension $3n$. As $(x_i,y_i)$ is a system of local coordinates of $TM$ in a neighbourhood
of $\pi_1(X)$, the triples $(\tm,\pi_1,TM)$ and $(\tm,\pi_2,M)$ have the structure of fibered
manifold. The triples  $(TT_2M,T(\pi_1),TTM)$ and $(TT_2M,T(\pi_2),TM)$ have also the 
structure of fibered manifold. We denote $\ttm=\tm-\{0\}$ and $\x$ the module of vector fields 
on $\tm$.

A tangent vector to $\tm$ is $\pi_1-$vertical (resp.$\pi_2-$vertical) if its projection by $T(\pi_1)$
(resp.$T(\pi_2)$) on $TTM$ (resp.$TM$) vanishes. If the elements of $T\tm$ are represented
locally by $(x,y,z,p,q,r)$, a $\pi_1-$vertical (resp.$\pi_2-$vertical) vector is represented locally 
by $(x,y,z,0,0,r)$(resp.$(x,y,z,0,q,r)$). We denote by $\vo$ (resp.$\vt$) the bundle of
$\pi_1-$vertical (resp.$\pi_2-$vertical) vectors. The vector bundles $\vo$ (resp.$\vt$) and
$\tm\times_MTM$ (resp.$\tm\times_{TM}TTM$), each of base $\tm$, are isomorphic

We have the following two exact sequences of vector bundles of base $\tm$:
$$0\longrightarrow\vt\stackrel{i_1}{\longrightarrow}T\tm\stackrel{s_1}{\longrightarrow}
\tm\times_MTM\longrightarrow0$$
$$0\longrightarrow\vo\stackrel{i_2}{\longrightarrow}T\tm\stackrel{s_2}{\longrightarrow}
\tm\times_{MT}TTM\longrightarrow0$$
where $i_1,s_1,i_2,s_2$ are defined locally by
\begin{eqnarray*}
i_1(x,y,z,p,q) &=&(x,y,z,0,p,q),\qquad s_1(x,y,z,p,q,r) = ((x,y,z),(x,p))\\
i_2(x,y,z,p)    &=&(x,y,z,0,0,p),\qquad s_2(x,y,z,p,q,r) =((x,y,z),(x,y,p,q))
\end{eqnarray*}
(These morphisms are defined globally also: for the details cf. [1]).

Let $h_1$ and $h_2$ be the inverses of the isomorphisms $\vo\longrightarrow\tm\times_MTM$
and $\vt\longrightarrow\tm\times_{TM}TTM$ [1] and define the vector $1$-forms:
$$J_1:T\tm\longrightarrow T\tm:J_1=i_2\circ h_1\circ s_1$$
$$J_2:T\tm\longrightarrow T\tm:J_2=i_1\circ h_2\circ s_2$$
$J_1$ and $J_2$ are called the canonical tensors of $\tm$. They play a very essential role
through the whole theory.

The canonical tensors $J_1$ and $J_2$ have the properties:
\begin{eqnarray}
& rank J_1=n,\ \  rank J_2=2n\nonumber\\
& Ker J_1=Im J_2=\vt,\ \  Ker J_2=Im J_1=\vo\\
& J_1^2=0,\ \ J_2^2=2J_1,\ \ J_1J_2=0=J_2J_1,\ \ J_1^3=0\\
& [J_1,J_1]=0,\ \ [J_2,J_2]=0,\ \ [J_1,J_2]=0
\end{eqnarray}

Let $\alpha$ be the section of the vector bundle $\tm\times_{TM}TTM$ (of base $\tm$) defined
by $\alpha=id_{\tm}\times_{TM}i$, where $i$ is the natural injection of $\tm$ in $TTM$
($i(J_0^2f)=J_0^1(J_0^1f))$, and let $\beta$ be the section of the vector bundle $\tm\times_MTM$
(of base $\tm$) defined by $\beta=id_{\tm}\times_M\pi_1$. With the aid of $\alpha$ and $\beta$
we define the vector fields on $\tm$:
$$C_1=i_2\circ h_1\circ\beta,\qquad C_2=i_1\circ h_2\circ\alpha$$
$C_1$ and $C_2$ are called the canonical vector fields of $\tm$. We have
\begin{equation}
J_1C_1=0,\ \ J_2C_1=0,\ \ J_1C_2=0,\ \  J_2C_2=2C_1,\ \ [C_1,C_2]=C_1
\end{equation}

We also have the following properties:
\begin{equation}
[C_1,J_1]=0,\ \ [C_1,J_2]=-J_1,\ \ [C_2,J_1]=-2J_1,\ \  [C_2,J_2]=-J_2
\end{equation}

The canonical vector field $C_2$ can be used to measure the homogeneity. A sclar $p$-form 
$\omega$, $p\ge0$, on $\tm$ is homogeneous of degree $r$ if $d_{C_2}\omega=r\omega$, 
where $d_{C_2}$ is the Lie derivative with respect to $C_2$. A vector $l$-form $L$, $l\ge0$,
on $\tm$ is homogeneous of degree $r$ if $[C_2,L]=(r-1)L$. We denote the homogeneity of
degree $r$ by the symbol $h(r)$.

A sclar $p$-form $\omega$, $p\ge1$, on $\tm$ is is $\pi_1$-semi-basic (resp. $\pi_2$-semi-
basic ) if $i_{J_1X}\omega=0$ (resp. $i_{J_2X}\omega=0$ for all $X\in\x$. A vector $l$-form 
$L$, $l\ge1$, on $\tm$ is is $\pi_1$-semi-basic (resp. $\pi_2$-semi-basic ) if $J_2L=0$ and 
$i_{J_1X}L=0$(resp.  $J_1L=0$ and $i_{J_2X}L=0$).

A semi-spray of type $1$ (resp. type $2$) is a vector field $S$ on $\tm$, differentiable on 
$\ttm$, such that $J_1S=C_1$ (resp. $J_2S=C_2$). A semi-spray $S$ of type $1$ (resp. 
type~$2$~) is a spray of type $1$ (resp. type $2$) if $S$ is homogeneous of degree $2$:
 $[C_2,S]=S$.

If $S$ is a semi-spray of type 1, we have
\begin{equation} 
J_1[J_1,S]=0,\ \ J_1[J_2,S]=J_1
\end{equation}

If $S$ is a semi-spray of type2, we have
\begin{equation} 
J_1[J_1,S]=0,\ \ J_1[J_2,S]=J_1,\ \ J_2[J_1,S]=2J_1,\ \ J_2[J_2,S]=J_2
\end{equation}

We also have
\begin{equation} 
[J_2,S]J_1=-2J_1\ \ [J_2,S]J_2=2[J_1,S]-J_2
\end{equation}

A (nonlinear) connection of type $1$ (resp. type $2$) on $\tm$ is a vector $1$-form $\G$ on
$\tm$, $C^\infty$ on $\ttm$, such that $J_1\G=J_1,\ \G J_2=-J_2$ (resp. $J_2\G=J_2,
\ \G J_1=-J_1$). $\G$ has the property that $\G^2=I$. A (nonlinear) connection $\G$ of type $1$ 
or of  type $2$ on $\tm$ is homogeneous if the vector $1$-form $\G$ is homogeneous of 
degree $1$: $[C_2,\G]=0$.

The horizontal and vertical projectors of $\G$ are $h=\frac{1}{2}(I+\G), \ v=\frac{1}{2}(I-\G)$
and we have the direct sum decomposition
$$T_z\tm=H_z(\tm)\oplus V_z(\tm) \ \  \forall z\in\tm,$$  
where, for a connection of type $1$
$$H_z(\tm)=T_{\pi_2(z)}M, \ \ V_z(\tm)=V^{\pi_2}_z(\tm)$$
and for a connection of type $2$
$$H_z(\tm)=T_{\pi_1(z)}TM, \ \ V_z(\tm)=V^{\pi_1}_z(\tm)$$
Hence a connection of type $1$ (resp. type $2$) is associated to the fibred manifold
$(\tm,\pi_2,M)$ (resp. $(\tm,\pi_1,TM)$).

For a connection of type $1$ on $\tm$, we have
\begin{equation}
J_1h=J_1, \ \ J_1v=0, \ \ hJ_2=0, \ \ vJ_2=J_2 
\end{equation}
and for a connection of type $2$ on $\tm$, we have
\begin{equation}
J_2h=J_2, \ \ J_2v=0, \ \ hJ_1=0, \ \ vJ_1=J_1 
\end{equation}

To each connection $\G$ of type $1$ (resp. type $2$) is canonically associated a semi-spray 
$S$ of type $1$ (resp. type $2$) given by $S=hS'$, where $S'$ is an arbitrary semi-spray of 
type $1$ (resp. type $2$). If $\G$ is homogeneous, $S$ is a spray (of the corresponding type).

To each spray $S$ of type $2$ are associated two homogeneous connections, one of type 
$1$, $\G_1$, and the other of type $2$, $\G_2$, wher,
\begin{eqnarray}
\G_1 &=& \frac{1}{3}\{2[J_2,S]+2[[J_1,S],S]-I\}\\
\G_2 &=& \frac{1}{3}\{2[J_2,S]+I\}
\end{eqnarray}
Fixing the given spray $S$ of type $2$, each one of the two connections $\G_1$ and 
$\G_2$ determines the other [1]; a reason for which we say that the two connections are
cojugates with respect to $S$.
\newpage

\Section{ Canonical Decomposition of Connections} 
\vspace{-0.5cm}
\begin{center}{\Large\bf of Type $1$ on $\tm$} 
\end{center}
\par
In [1], CATZ proved the so-called theorem of canonical decomposition of connections of type
$2$ on $\tm$, namely, "Given a spray $S$ of type $2$ on $\tm$ and a vector $1$-form $T$
on $\tm$, $\pi_1$-semi-basic such that $T(S)=0$, there exists a unique homogeneous 
connection $\G$ of type $2$ on $\tm$ whose associated spray of type $2$ is $S$ and whose 
strong torsion is $T$. It is given by $\G=\frac{1}{3}\{2[J_2,S]+T+I\}$".

CATZ succeeded to prove this result by making use of the many interesting identitis 
satisfied by sprays of type $2$; cf. for example equations (7) and (8). But for sprays of type 
$1$ we have very few identities, and so such a decomposition is not evident. Nevertheless, 
we try in this section to establish some sort of decomposition theorem for connections of 
type $1$ on $\tm$.
\begin{definition}
Let $\G$ be a connection of type $1$ on $\tm$. The weak torsion of $\G$ is the $\pi_2$-semi
-basic vector $2$-form on $\tm$ given by \ $t=[J_1,\G]$.
\par The strong torsion of $\G$ is the vector $1$-form on $\tm$ given by:
$$T=i_St-[C_2,\G],$$
where $S$ is an arbitrary semi-spray of type $1$ on $\tm$. If $\G$ is homogeneous, then 
$T=i_St$. It can be shown that the definition of strong torsion is independent of the choice 
of $S$.
\end{definition}
\begin{proposition}
The strong torsion $T$ of a connection $\G$ of type $1$ on $\tm$ can be expressed in the
form: 
\begin{equation}
T=-J_2v+2[S,J_1]+[C_1,\G]-[C_2,\G],
\end{equation}
where $S$ is the semi-spray of type $1$ associated to $\G$.
\par The strong torsion $T$ has the following properties:\newline
(a) $T$ is $\pi_2$-semi-basic,\newline
(b) $T(S)=-2v[C_2,S]$ and $T(S)=0$ if $\G$ is homogeneous.
\end{proposition}
{\bf Proof}. We have
\begin{eqnarray*}
T &=& i_St-[C_2,\G]=i_S[J_1,\G]-[C_2,\G]\\
   &=& J_1[\G,S]+\G[J_1,S]-[J_1,\G S]-[\G,J_1S]-[C_2,\G]
\end{eqnarray*}
Taking for $S$ the semi-spray of type $1$ associated to $\G$: $S=hS'$, where $S'$ is an
arbitrary semi-spray of type $1$. Thus $\G S=\G hS'=hS'=S$. From which we have
\begin{equation}
T=J_1[\G,S]+\G[J_1,S]-[J_1,S]+[C_1,\G]-[C_2,\G]
\end{equation}
We calculate the term $J_1[\G,S]$:\newline
As $J_1(\G-I)=0$,\ $\G-I\in Ker J_1$. Thus, by (1), $\G-I\in Im J_2$ and $\G-I$ can be written as 
$\G-I=J_2L$, where $L$ is a vector $1$-form on $\tm$. Then
$$J_1[\G,S]=J_1[\G-I,S]=J_1[J_2L,S]=J_1\{[J_2,S]L+J_2[L,S]\}$$
Using (2) and (6), we get
\begin{equation}
J_1[\G,S]=J_1[J_2,S]L=J_1L=\frac{1}{2}J^2_2L=\frac{1}{2}J_2(\G-I)=J_1[\G,S]=-J_2v
\end{equation}
We now calculate the term $\G[J_1,S]$:\newline
As $J_1[J_1,S]=0$ by (6), we have $[J_1,S]\in Ker J_1=Im J_2$, and since $\G J_2=-J_2$,
we get
\begin{equation}
\G [J_1,S]=-[J_1,S]
\end{equation}
Substituting from (15) and (16) into (14), we get the required.

(a) To prove that $T$ is $\pi_2$-semi-basic, we have to show that $\ J_1T=0$ and $\ T(J_2X)
=0$ \ for all $X\in\x$.
\begin{eqnarray*}
J_1TX &=& -J_1J_2vX+2J_1[S,J_1]X+J_1[C_1,\G]X-J_1[C_2,\G]X\\
           &=& J_1[C_1,\G]X-J_1[C_2,\G]X,\mbox{  by (2) and (6)}\\
           &=& J_1[C_1,\G X]-J_1\G[C_1,X]-J_1[C_2,\G X]+J_1\G[C_2,X]\\
           &=& J_1[C_1,\G X]-J_1[C_1,X]-J_1[C_2,\G X]+J_1[C_2,X]
\end{eqnarray*}
Now, using (5), we get
\begin{eqnarray*}
         0 &=& [C_1,J_1]\G X=[C_1,J_1\G X]-J_1[C_1,\G X]\Longrightarrow J_1[C_1,\G X]
                   =[C_1,J_1 X],\\
-2J_1X &=&[C_2,J_1]\G X=[C_2,J_1\G X]-J_1[C_2,\G X]\Longrightarrow J_1[C_2,\G X]
                   =[C_2,J_1X]+2J_1X
\end{eqnarray*}
Substituting the above two equations in the expression of $J_1TX$, we obtain
\begin{eqnarray*}
J_1TX &=& [C_1,J_1X]-J_1[C_1,X]-[C_2,J_1 X]+J_1[C_2,X]-2J_1X\\
           &=& [C_1,J_1] X-[C_2,J_1]X-2J_1X]=0, \mbox{  again by (5)}
\end{eqnarray*}
Hence $J_1T=0$. On the other hand,
$$T(J_2X)=-J_2v(J_2X)+2[S,J_1]J_2X+[C_1,\G]J_2X-[C_2,\G]J_2X$$
Using (9), (2) and (6), we get after some calculation
$$J_2v(J_2X)=2J_1X,\quad [S,J_1]J_2X=J_1X,\quad [C_1,\G]J_2X=-2h[C_1,J_2X]$$
From (5), we have
$$0=[C_1,J_1]J_2X=[C_1,J_1J_2 X]-J_1[C_1,J_2X]=-J_1[C_1,J_2X],$$
from which $[C_1,J_2X]\in Ker J_1=Im J_2=\vt$, consequently, $h[C_1,J_2X]=0$. Hence
$[C_1,\G]J_2X]=0$. Similarly, $[C_2,\G]J_2X=0$. We thus coclude that $T(J_2X)=0$.

(b)$\ \ \ T(S)=-J_2vS+2[S,J_1]S+[C_1,\G]S-[C_2,\G]S.\ $
Using the fact that $\G S=S$ and that $J_1S=C_1$, we get
$$J_2vS=0,\quad[S,J_1]S=-[C_1,S],\quad[C_1,\G]S=[C_1,S]-\G[C_1,S],\quad[C_2,\G]S
=2v[C_2,S]$$
Hence, we have
$$T(S)=-[C_1,S]-\G[C_1,S]-2v[C_2,S]=-2h[C_1,S]-2v[C_2,S]$$
From the identity that $[C_1,J_1]S=0$, we deduce that $[C_1,S]$ is $\pi_2$-vertical, and
hence $h[C_1,S]=0$. Thus $T(S)=-2v[C_2,S]$.
\par If $\G$ is homogeneous, the semi-spray $S$ associated to $\G$ is a spray:$[C_2,S]=S$,
and in this case $T(S)=0$.\quad$\Box$

\begin{remark}
{\em As $v=\frac{1}{2}(I-\G)$, the formula (13) can be written as
$$T=\frac{1}{2}\{J_2,\G-J_2-4[J_1,S]+2[C_1,\G]-2[C_2,\G]\}$$
and if $\G$ is homogeneous such that $[C_1,\G]=0$, this formula takes the form
\begin{equation}
T=\frac{1}{2}\{J_2\G-J_2-4[J_1,S]\}
\end{equation}
}
\end{remark}

\par We have seen that to a homogeneous connection $\G$ of type $1$ on $\tm$ there are
canonically associated:\newline
--- a spry of type $1$ on $\tm$, $\ S=hS'$ ($S'$: arbitrary spray of type $1$),\newline
--- a vector $1$-form $T$ on $\tm$, $\pi_2$-semi-basic, called strong torsion of $\G$, such
that $T(S)=0$. 

We are now in a position to announce a canonical decomposition version for connections
of type $1$ on $\tm$ (cf. Eq. (17)):

\begin{theorem}
Given a spray $S$ of type $1$ on $\tm$ and a vector $1$-form $T$ on $\tm$, $\pi_2$-semi-basic
such that $T(S)=0$, there exists a homogeneous connection $\G$ of type $1$ on $\tm$ whose
associated spray of type $1$ is $S$, whose strong torsion is $T$ and having the property
that $[C_1,\G]=0$. This connection is determined by the formula:
\begin{equation}
J_2\G=2T+J_2+4[J_1,S]
\end{equation}
\end{theorem}

\begin{remark}
{\em The expression (18) may be written as 
$J_2(\G-I)=2(T+2[J_1,S])$.
Since $J_1(T+2[J_1,S])=0$, then $T+2[J_1,S]\in Ker J_1=Im J_2$; so there exists a vector 
$1$- form $L$ on $\tm$ such that $T+2[J_1,S]=J_2L$. From which 
$J_2(\G-I-2J_2L)=0$.
Thus there exists another vector $1$-form $K$ on $\tm$such that $\G-I-2J_2L=J_1K$. Hence
$$\G=I+J_1K+J_2L,$$
where $(J_1K+2J_2L)J_2=-2J_2$.
}
\end{remark}

\Section{ Nonlinear Connections Induced by a } 
\vspace{-0.5cm}
\begin{center}{\Large\bf Linear Connectionon $\tm$} 
\end{center}
\par
In this section we show that a linear connection on $\tm$ satisfying certain conditions
induces nonlinear connections of type $1$ and of of type $2$ on $\tm$ and we investigate 
these induced connections.

We denote by $\xo$ (resp. $\xt$) the set of $\pi_1$-vertical (resp. $\pi_2$-vertical) vector
fields on $\tm$.

Let $D$ be a linear connection on $\tm$. Consider the mappings
$$DC_1: \x\longrightarrow\x:X\longmapsto D_XC_1$$
$$DC_2: \x\longrightarrow\x:X\longmapsto D_XC_2$$
We denote $\phi=DC_1\mid_{\xo}$ and $\psi=DC_2\mid_{\xt}$. 
IF $DJ_1=0$, then $Im\, DC_1\subseteq\xo$ and $\phi$ is thus a mapping of $\xo$ into itself.
IF $DJ_2=0$, then $Im\, DC_2\subseteq\xt$ and $\psi$ is thus a mapping of $\xt$ into itself.
We have 
$$DC_1\circ J_1=\phi\circ J_1,\quad DC_2\circ J_2=\psi\circ J_2$$
It should be noticed that $DJ_2=0\Longrightarrow DJ_1=0$, but the converse is not true. 
\begin{definition}
A linear connection $D$ on $\tm$ is said to be $J_1$-regular if:\newline
(a) $DJ_1=0$\newline
(b) the mapping $\phi:\xo\longrightarrow\xo:X\longmapsto D_XC_1$ is an isomorphism.
\par
A linear connection $D$ on $\tm$ is said to be $J_2$-regular if:\newline
(a) $DJ_2=0$\newline
(b) the mapping $\psi:\xt\longrightarrow\xt:X\longmapsto D_XC_2$ is an isomorphism.
\end{definition}

The followig result shows that a $J_1$-regular ( resp. $J_2$-regular) linear connection on
$\tm$ induces a nonlinear connection of type $2$ (resp. type $1$) on $\tm$. This generalizes
a known result [4] (related to $TM$) to the tangent bundle of order $2$.
\begin{theorem} 
Let $D$ be a linear connection on $\tm$.\newline
(a) If $D$ is $J_1$-regular, the vector $1$-form on $\tm$:
$$\G_2=I-2\phi^{-1}\circ DC_1$$
is a connection of type $2$ on $\tm$.\newline
(b) If $D$ is $J_2$-regular, the vector $1$-form on $\tm$:
$$\G_1=I-2\psi^{-1}\circ DC_2$$
is a connection of type $1$ on $\tm$.
\par The nonlinear connection $\G_2$(resp. $\G_1$) on $\tm$ is said to be induced by the
linear connection $D$ on $\tm$. 
\end{theorem} 
{\bf Proof.}
\par (a) We have $\ \ \ J_2\G_2=J_2-2J_2(\phi^{-1}\circ DC_1)$\newline
By the hypothesis, $(\phi^{-1}\circ DC_1)X\in\xo$ for all $X\in\x$, from which 
$J_2(\phi^{-1}\circ DC_1)=0$. Thus we have $J_2\G_2=J_2$. On the other hand
$$\G_2J_1=J_1-2(\phi^{-1}\circ DC_1)J_1$$
Now, $(\phi^{-1}\circ DC_1)\circ J_1=\phi^{-1}\circ( DC_1\circ J_1)=\phi^{-1}\circ(\phi\circ J_1)
=J_1$, and hence $\G_2J_1=-J_1$.
\par (b) can be proved analogously.\quad$\Box$
\begin{proposition}~~\newline
(a) The connection $\G_2$ of type $2$ induced by a $J_1$-regular linear connection $D$ is
homogeneous if, and only if, $\ [C_2,DC_1]=0$\newline
(b) The connection $\G_1$ of type $1$ induced by a $J_2$-regular linear connection $D$ is
homogeneous if, and only if, $\ [C_2,DC_2]=0$
\end{proposition}
{\bf Proof.}
\par (a)Let $D$ be a $J_1$-regular linear connection on $\tm$ such that $\ [C_2,DC_1]=0$.
We have:
\begin{eqnarray*}
\frac{1}{2}[C_2,\G_2] &=& -[C_2,\phi^{-1}\circ DC_1]
                                         =-[C_2,\phi^{-1}]\circ DC_1-\phi^{-1}\circ [C_2,DC_1]\\
                                  &=& -[C_2,\phi^{-1}]\circ DC_1,
\end{eqnarray*}
On the other hand,
\begin{eqnarray*}
0 &=& \phi^{-1}[C_2, DC_1](\phi^{-1}\circ DC_1)X\quad\forall X\in\x\\
   &=& \phi{^-1}[C_2,(DC_1\circ\phi^{-1}\circ DC_1)X]-(\phi^{-1}\circ DC_1)[C_2,(\phi^{-1}
          \circ DC_1)X]\\
\end{eqnarray*}
But since $\ DC_1\circ\phi^{-1}\mid_{\xo}=\phi^{-1}\circ DC_1\mid_{\xo}=Id_{\xo}$, then
\begin{eqnarray*}
0 &=& \phi^{-1}[C_2,(DC_1)X]-[C_2,(\phi^{-1}\circ DC_1)X]\\
   &=& \phi^{-1}[C_2,(DC_1)X]-[C_2,\phi^{-1}](DC_1)X-\phi^{-1}[C_2,(DC_1)X]\\
   &=& -[C_2,\phi^{-1}](DC_1)X\\
   &=& \frac{1}{2}[C_2,\G_2]X
\end{eqnarray*}
The converse is proved analogously.
\par (b) can be proved using the same arguments.\quad$\Box$
~~\newline
\par
It is proved in [1] that any linear connection $D$ on $\tm$ with no torsion and such that
$DJ_2=0$ is $J_2$-regular and thus induces a connection of type $1$ on $\tm$. This result 
fails to be true for $J_1$-regular connections, and instead we have
\begin{proposition}
On $\tm$, there are no $J_1$-regular linear connections with no torsion.
\end{proposition}
{\bf Proof.}
Let $D$ be a $J_1$-regular linear connection on $\tm$ without torsion. We show that the 
hypothesis that $T=0$ and $DJ_1=0$ imply that $\phi$ can not be an isomorphism. For all 
$X, Y\in\x$, we have
\begin{eqnarray*}
D_XJ_1Y &=& J_1D_XY\\
      D_XY &=& D_YX+[X,Y]
\end{eqnarray*}
Then,
\begin{eqnarray*}
D_{J_1X}C_1 &=& D_{C_1}J_1X+[J_1X,C_1]=J_1D_{C_1}X+[J_1X,C_1]\\
                      &=& J_1D_{C_1}X-[C_1,J_1]X+J_1[X,C_1]\\
                      &=& J_1(D_{C_1}X+[X,C_1]), \quad\mbox{ by (5)}\\
                      &=& J_1D_XC_1=D_XJ_1C_1=0,\quad\mbox{ by (4)}
\end{eqnarray*}
From which, $\ \phi(J_1X)=D_{J_1X}C_1=0\quad\forall X\in\x\ $, which means that $\phi=0$
\quad$\Box$

\par
Let $D$ be a $J_1$-regular linear connection on $\tm$ such that $[C_2,DC_1]=0$. The 
nonlinear connection
$$\G_2=I-2\phi^{-1}\circ DC_1$$
of type $2$ induced by $D$ is thus homogeneous. Let $S$ be the spray of type $2$ 
associated to $\G_2$. According to (12), this spray $S$ of type $2$ defines another 
connection 
\begin{equation}
\o{\G}_2=\frac{1}{3}\{2[J_2,S]+I\}
\end{equation}
of type $2$ on $\x$, with no strong torsion and whose associated spray of type $2$ is $S$
itself. We conclude that the spray $S$ of type $2$ is simultaneously associated to both 
connections $\G_2$ and $\o{\G}_2$: $\ h_2S=\o{h}_2S=S$ or $\ \G_2S=\o{\G}_2S=S$, where 
$\ h_2$ (resp. $\o{h}_2$) is the horizontal projector of $\G_2$ (resp.$\o{\G}_2$).

The strong torsion of $\G_2$ is given by
$$T=3\G_2-I-2[J_2,S]$$
But from (19), we have
$$2[J_2,S]=3\o{\G}_2-I$$
Hence,
$$T=3(\G_2-\o{\G}_2),$$
from which the following result.
\begin{proposition}
Let $D$ be a $J_1$-regular linear connection on $\tm$ such that $[C_2,DC_1]\\=0$. Let $\G_2$
be th econnection of  type $2$ induced by $D$ and $\o{\G}_2$ the connection of type $2$ 
defined by the spray $S$ of type $2$ associated to $\G_2$. The two connections $\G_2$ and 
$\o{\G}_2$ coincide if, and only if, $\G_2$ has no strong torsion. In this case we have 
$\ DC_1=\frac{1}{3}\phi\circ(I-[J_2,S])$.
\end{proposition}
\begin{remark}
{\bf The spray $S$ of type $2$ associated to $\G_2$ (and to $\o{\G}_2)$ defines also on $\tm$ 
a homogeneous connection $\o{\G}_1$  of type $1$ given by (cf. Eq. (11)):
$$\o{\G}_1=\frac{1}{3}\{2[J_2,S]+2[[J_1,S],S]-I\}$$ 
The two connections $\o{\G}_1$ and $\o{\G}_2$ are cojugates with respect to $S$.
}
\end{remark}

\Section{ Finslerian Forms, Sprays and} 
\vspace{-0.5cm}
\begin{center}{\Large\bf Connectionon $\tm$} 
\end{center}
\par 
In this section we define and investigate the Finslerian forms $\O$ on $\tm$. We also
prove the existence and uniqueness of a spray of typ $2$, a connection of type $2$ and a
connection of type $1$, canonically constructed on $(\tm,\O)$.
\begin{lemma}
For all semi-spray $S$ of type $2$ on $\tm$, we have
$$[i_S,d_{J_2}]=d_{C_2}-i_{[S,{J_2}]}$$
\end{lemma}
{\bf Proof.} See [2], page 352.
\begin{proposition}
Let $\omega$ be a scalar $p$-form on $\tm$, $\pi_1$-semi-basic \em{(}resp. $\pi_2$
-semi-basic\em{)}, homogeneous of degree $r$. For all semi-spray $S$ of type $2$ on $\tm$,
we have
$$[i_S,d_{J_2}]\omega=(r+p)\omega$$
or
$$i_Sd_{J_2}\omega+d_{J_2}i_S\omega=(r+p)\omega$$
\end{proposition}
{\bf Proof.} We have from the hypothesis
$$d_{C_2}\omega=r\omega,\quad i_{J_1Z}\omega=0 (\mbox{resp.\  }i_{J_2Z}\omega=0)
\quad\forall Z\in\x.$$
After the precedent lemma, we have
$$i_Sd_{J_2}\omega+d_{J_2}i_S\omega=d_{C_2}\omega-i_{[S,J_2]}\omega=r\omega
-i_{[S,J_2]}\omega$$
It remains to show that  $i_{[S,J_2]}\omega=-p\omega$.
\begin{equation}
(i_{[S,J_2]}\omega)(Z_1,\cdots,Z_p)=\sum_{i=1}^{p}\omega(Z_1,\cdots,[S,J_2]Z_i,\cdots,Z_p)
\end{equation}
\par Case 1. $\omega$ is $\pi_1$-semi-basic: $i_{J_1Z}\omega=0$.\newline
We have from (7), $\ J_2[S,J_2]Z_i=-J_2Z_i$, or
$$J_2([S,J_2]Z_i+Z_i)=0$$
Using (1), we can write
$$[S,J_2]Z_i+Z_i=J_1W,  \mbox{ for some }W\in\x$$
Substituting in (20), we get
\begin{eqnarray*}
(i_{[S,J_2]}\omega)(Z_1,\cdots,Z_p) &=& \sum_{i=1}^{p}\omega(Z_1,\cdots,J_1W-Z_i,
                                                                  \cdots,Z_p)\\
                                                          &=& -\sum_{i=1}^{p}\omega(Z_1,\cdots,Z_i,\cdots;Z_p)\\
                                                          &=& -p\omega(Z_1,\cdots,Z_p)
\end{eqnarray*}
\par Case 2. $\omega$ is $\pi_2$-semi-basic: $i_{J_2Z}\omega=0$.\newline
Using (6) and (1), one can write
$$[S,J_2]Z_i+Z_i=J_2U, \mbox{ for some }U\in\x.$$
Substituting in (20), we get $i_{[S,J_2]}\omega=-p\omega$.
\par Hence the result. \quad$\Box$\newline
\par In the course of the above proof we have proved the
\begin{corollary}
For all scalar $p$-form $\omega$ on $\tm$, $\pi_1$-semi-basic \em{(}resp. $\pi_2$-semi-
basic\em{)}, we have
$$i_{[S,J_2]}\omega=-p\omega,$$
where $S$ is any arbitrary semi-spray of type $2$ on $\tm$.
\end{corollary}
\par Proposition 5 implies directly the following result.
\begin{theorem}
Any $p$-form $\omega$ on $\tm$, $\pi_1$-semi-besic \em{(}resp. $\pi_2$-semi-basic\em{)}, 
homogeneous of degree $r\not=-p$, which is $d_{J_2}$-closed is $d_{J_2}$-exact and
$$\omega=\frac{1}{r+p}\, d_{J_2}i_S\omega,$$
where $S$ is any arbitrary semi-spray of type $2$ on $\tm$.
\end{theorem}
~~\newline
\par
Let $\O$ be a scalar $2$-form of maximal rank $3n$ on $\tm$. It is clear that the formula
$i_Z\O=\omega$ defines an isomorphism between the vector space of vector fields $Z$ on
$\tm$ and the vector space of scalar $1$-forms $\omega$ on $\tm$. If, moreover, $\O$ is
homogeneous of degree $1$ ($d_{C_2}\O=\O$), this formula defines an isomorphism between
the vector fields $h(2)$ and the $1$-forms $h(2)$. In fact,
$$i_{[C_2,Z]}\O=-[i_Z,d_{C_2}]\O=d_{C_2}i_Z\O-i_Zd_{C_2}\O=d_{C_2}i_Z\O-i_Z\O
=d_{C_2}\omega-\omega$$
Thus, we have
$$\omega:h(2)\Longleftrightarrow d_{C_2}\omega=2\omega\Longleftrightarrow i_{[C_2,Z]}\O
=\omega=-i_Z\O\Longleftrightarrow[C_2,Z]=Z\Longleftrightarrow Z:h(2)$$
\par We shall now consider certain type of these $2$-forms.
\begin{definition}
A scalar $2$-form $\O$ of maximal rank $3n$ on $\tm$, homogeneous of degree $1$, such
that $i_{J_2}\O=0$ will be called a Finslerian form on $\tm$.
\end{definition}
\par The following result can be proved easily.
\begin{proposition}
Any Finslerian form $\O$ on $\tm$ induces a metric $g$ on $\vt$ given by
$$g(J_2X,J_2Y)=\O(J_2X,Y)\qquad\forall X,Y\in\x.$$
such a metric will be called Finslerian metric on $\tm$ induced by $\O$.
\end{proposition}
\begin{proposition}
Let $\O$  be a Finslerian form on $\tm$ and let $\omega=i_{C_2}\O$. To each $1$-form 
$\alpha$ on $\tm$, homogeneous of degree $2$ and verifying $J_2\alpha=-\omega$, there 
correspods a spray $S$ of type $2$ defined by $i_S\O=\alpha$.
\end{proposition}
{\bf Proof.} We have $\ \alpha=i_S\O,\quad\omega=i_{C_2}\O,\quad\omega=-J_2\alpha$.
To prove that $S$ is a spray of type $2$, we have to show that $J_2S=C_2$ and $[C_2,S]=S$.
$$i_{J_2S}\O=[i_S,i_{J_2}]\O=i_Si_{J_2}\O-i_{J_2}i_S\O=-i_{J_2}i_S\O=-i_{J_2}\alpha
=-J_2\alpha=\omega=i_{C_2}\O,$$
from which $J_2S=C_2$. On the other hand, since $\O$ is $h(1)$ and since $\alpha$ is 
$h(2)$, then $S$ is $h(2)$: $[C_2,S]=S$.
\quad$\Box$
~~\newline
\par
Now,we define a function $E$ on $\tm$ by
$E=\frac{1}{2}\O(C_2,S)=\frac{1}{2}i_Si_{C_2}\O,$
where $S$ is any arbitrary spray of type $2$ on $\tm$. This definition does not depend on 
the choice of $S$. In fact, let $S'$ be another spray of type $2$ on $\tm$. Then, $J_2S'=J_2S
=C_2$, from which $J_2(S'-S)=0$. Thus we have
$$0=(i_{J_2}\O)(S,S'-S)=\O(J_2S,S'-S)+\O(S,J_2(S'-S))=\O(C_2,S'-S),$$
from which $\O(C_2,S')=\O(C_2,S)$ and the function $E$ is thus well-defined.
\par We show that $E$ is homogeneous of degree $2$: $d_{C_2}E=C_2.E=2E$. Since 
$\O$ is $h(1)$ , we have 
\begin{eqnarray*}
2E &=& \O(C_2,S)=(d_{C_2}\O)(C_2,S)=C_2.\O(C_2,S)-\O([C_2,C_2],S)-\O(C_2,[C_2,S])\\
     &=& d_{C_2}\O(C_2,S)-\O(C_2,S)=2d_{C_2}E-2E,
\end{eqnarray*}
from which $d_{C_2}E=2E$ and $E$ is $h(2)$. We can thus set the
\begin{definition}
Let $\O$ be a Finslerian form on $\tm$. We call energy function associated to $\O$ the 
function $E$, homogeneous of degree $2$, on $\tm$ defined by
$$E=\frac{1}{2}\O(C_2,S)=\frac{1}{2}i_Si_{C_2}\O,$$
where $S$ is any arbitrary spray of type $2$ on $\tm$.
\end{definition}
\begin{remark}
Taking Proposition 6 into account, the energy function $E$ may be defined by
$$E=\frac{1}{2}\o{g}(C_2,C_2).$$  
\end{remark}
\par The following two theorems generalize the main results of [6] to the tangent bundle of
orde $2$.
\begin{theorem}
Let $\O$ be a Finslerian form on $\tm$. If the $1$-form $i_{C_2}\O$ on $\tm$ is $d_{J_2}$-closed,
then $i_{C_2}\O=d_{J_2}E$ and the vector field $G$ on $\tm$ defined by $i_G\O=-dE$ is a 
spray of type $2$ on $\tm$, called the canonical spray of $(\tm,\O)$.
\end{theorem}
{\bf Proof.}
We show firstly that the form $i_{C_2}\O$ is $(h)$.
$$(d_{C_2}i_{C_2}\O)(Y)=C_2.\O(C_2,Y)-\O(C_2,[C_2,Y])=(d_{C_2}\O)(C_2,Y)=\O(C_2,Y)
=(i_{C_2}\O)Y,$$
since $\O$ is $h(1)$. Hence $d_{C_2}(i_{C_2}\O)=i_{C_2}\O$ and the form $i_{C_2}\O$ 
is $h(1)$.
\par Next, we show that $i_{C_2}\O$ is $\pi_1$-semi-basic.
$$0=(i_{J_2}\O)(S,J_1X)=\O(J_2S,J_1X)=\O(C_2,J_1X)=(i_{C_2}\O)(J_1X).$$
\par Now, $i_{C_2}\O$ is a $1$-form on $\tm$, $\pi_1$-semi-basic, homogeneous of degree 1
which is, by hypothesis, $d_{J_2}$-closed. We can thus apply Theorem 3 to the form 
$i_{C_2}\O$ to obtain
$$i_{C_2}\O=\frac{1}{2}\,d_{J_2}i_{C_2}\O=\frac{1}{2}\,d_{J_2}\O(C_2,S)=d_{J_2}E.$$
Hence, $i_{C_2}\O=d_{J_2}E=i_{J_2}dE=-i_{J_2}(-dE)$. Setting $\omega=i_{C_2}\O$ and
$\alpha=-dE$, we get $\omega=-J_2\alpha$, with $\alpha$ homogeneous of degree $2$ as
one can easily show. By Proposition~~7, to the form $\alpha$ on $\tm$ there corresponds a 
spray $S$ of type $2$ on $\tm$ defined by $i_G\O=\alpha=-dE$, which completes the proof.
\quad$\Box$
\begin{theorem}
A Finslerian form $\O$ on $\tm$ of energy $E$ can be written in the form
$$\O=dd_{J_2}E+\Theta,\quad\mbox{where } \Theta=i_{C_2}d\O.$$
If, moreover, $\O$ is closed, it is exact and
$$\O=dd_{J_2}E.$$
\end{theorem}
{\bf Proof.}
Since $\O$ is homogeneous of degree $1$,
$$\O=d_{C_2}\O=di_{C_2}\O+i_{C_2}d\O=dd_{J_2}E+i_{C_2}d\O.\quad\Box$$
\begin{proposition}
The canonical spray $G$ on $\tm$ has the properties:\newline
\em{(}a\em{)} $d_GE=0,$\newline
\em{(}b\em{)} b) $d_G\O=0$, if $\O$ is closed.
\end{proposition}
{\bf Proof.}
The fact that $i_G\O=-dE$ implies that $\,0=-i_GdE=-d_GE\,$ and that 
$\,d_G\O=i_Gd\O+di_G\O=di_G\O=-d^2E=0.$
\quad$\Box$
\begin{remark}
{\em We can prove that the $1$-form $i_{C_1}\O$ on $\tm$ is homogeneous of degree zero, 
$\pi_1$-semi-basic and $\pi_2$-semi-basic. Thus, according to Theorem 3, if $i_{C_1}\O$ 
is $d_{J_2}$-closed, we have $\,i_{C_1}\O=d_{J_2}i_Si_{C_1}\O=d_{J_2}\O(C_1,S)\,,$
where $S$ is any arbitrary spray of type $2$ on $\tm$. Now by (4), 
$\,0=(i_{J_2}\O)(C_2,S)=\O(J_2C_2,S)+\O(C_2,J_2S)=2\O(C_1,S)\,$.
Hence, if $i_{C_1}\O$ is $d_{J_2}$-closed, then $i_{C_1}\O=0$. We conclude that the 
reasoning employed for $i_{C_2}\O$ can not be applied to $i_{C_1}\O$.
}
\end{remark}
\par
Now, the theorem of canonical decomposition of connections of type $2$ (cf. $\S$ 2), for the case
$S=G$ (canonical spray of type $2$ of $(\tm,\O)$) and $T=0$, gives us the following result, 
which generalizes the fundamental theorem of Finsler geometry [3] to the tangent bundle 
of order $2$.
\begin{theorem}
On $(\tm,\O)$, there exists a unique homogeneous connection $\G_2$ of type $2$ on $\tm$ 
with no strong torsion and whose associated spray $G$ of type $2$ is defined by 
$i_G\O=-dE$. This connection is given by
\begin{equation}
\G_2=\frac{1}{3}\,\{2[J_2,G]+I\}
\end{equation}
and is called the canonical connection of type $2$ of $(\tm,\O)$.
\end{theorem}
As the canonical spray $G$ of $(\tm,\O)$ is a spray of type $2$, it defines another connection
of type $1$ on $(\tm,\O)$ according to (11). Hence we have
\begin{theorem}
Let $G$ be the canonical spray of $(\tm,\O)$.
On $(\tm,\O)$, there exists a unique homogeneous connection $\G_1$ of typi $1$ on $\tm$ 
associated to $G$. This connection is given by
$$\G_1=\frac{1}{3}\,\{2[J_2,G]+2[[J_1,G],G]-I\}$$
and is called the canonical connection of type $1$ on $(\tm,\O)$.
\par The two connections $\G_1$ and $\G_2$ are conjugates with respect to $G$.
\end{theorem}
{\bf Proof.} The existence of $\G_1$ (homogeneous) is insured by Eq (11) and the uniqueness 
of $\G_1$ is insured by the  uniqueness of $G$. In fact, if $G'$ is another spray of type $2$ 
such that $i_{G'}\O=-dE$, then $i_{G'-G}\O=0$ and $G'=G$, since $\O$ is of maximal rank. 
The fact that the two connections $\G_1$ and $\G_2$ ar defined by the same spray $G$ of 
type $2$ (according to (11) and (12)) shows that they are conjugates with respect to $G$.
\quad$\Box$
\par
~~\newline
{\bf Concluding Remark.} Comparing Theorem 6 and Proposition 4, we find that the 
connection $\G_2=\o{\G}_2$ of Proposition 4 (Eq. (19)) and the connection $\G_2$ of 
Theorem~~6 (Eq. (21)) have the same form (but with different associated sprays) and enjoy
similar properties. The spray $S$ of type $2$ in (19) is defined via a $J_1$-regular linear 
connection $D$ on $\tm$ verifying $[C_2,DC_1]=0$ while the spray $G$ of type $2$ in (21) 
is defined via  a Finslerian form $\O$ on $\tm$. The two geometric structures $(\tm,D)$ and  
$(\tm,\O)$ are thus strongly related.
\par It seems to us that endowing $\tm$ with a $J_1$-regular linear connection $D$ verifying
$[C_2,DC_1]=0$ is equivalent in some sense to endowing $\tm$ with a Finslerian form $\O$.
This problem may be the object of a forthcoming work.


\end{document}